\documentclass{ifacconf}

\usepackage{amsmath,amssymb}
\usepackage{mathtools}
\usepackage{tikz}
\usepackage{tikz-3dplot}
\usetikzlibrary{arrows.meta,calc,positioning}
\usepackage[normalem]{ulem}

\usepackage{graphicx}      
\usepackage{subcaption}    
\usepackage{float}         
\usepackage{natbib}        
\usepackage{acronym}       
\usepackage{todonotes}
\usepackage{enumitem}      
\acrodef{pbh}[PBH]{Popov-Belevitch-Hautus} 
\acrodef{uavs}[UAVs]{unmanned aerial vehicles}
\acrodef{mass}[MASs]{multi-agent systems}
\acrodef{lti}[LTI]{linear time-invariant}
\acrodef{uav}[UAV]{unmanned aerial vehicle}
\acrodef{idr}[IDR]{\emph{infinitesimal distance rigidity}}
\acrodef{rbm}[RBM]{\emph{rigid body motion}}
\acrodef{rbms}[RBMs]{\emph{rigid body motions}}
\acrodef{fvt}[FVT]{\emph{final value theorem}}
\acrodef{svd}[SVD]{\emph{singular value decomposition}} 
\acrodef{evd}[EVD]{\emph{eigen value decomposition}} 
\acrodef{siso}[SISO]{\emph{single-input single-output}}
\acrodef{mimo}[MIMO]{\emph{multi-input multi-output}} 
\acrodef{odes}[ODEs]{\emph{ordinary differential equations}} 
\acrodef{psd}[PSD]{\emph{positive semi-definite}} 
\acrodef{nsd}[NSD]{\emph{negative semi-definite}} 

\newcommand{\Gcal}{{\mathcal{G}}}
\newcommand{\Vcal}{{\mathcal{V}}}
\newcommand{\Ecal}{{\mathcal{E}}}

\newcommand{\rank}{{\mathrm{rank}}}
\newcommand{\image}{{\mathrm{Im}}}

\begin{document}
\raggedbottom
\begin{frontmatter}

\title{The Geometry of Hidden Modes in Distance-Based Formation Control\thanksref{footnoteinfo}} 

\thanks[footnoteinfo]{This work was supported by the Israel Science Foundation grant no. 453/24 and the Gordon Center for Systems Engineering.}

\author[First]{Solomon Goldgraber Casspi} 
\author[First]{Daniel Zelazo}

\address[First]{Stephen B. Klein Faculty of Aerospace Engineering, Technion - Israel Institute of Technology, Haifa, Israel (e-mail: solomon.g@campus.technion.ac.il, dzelazo@technion.ac.il).}

\begin{abstract}
This paper presents a geometric input-output analysis of hidden modes in distance-based
formation control. We study the linearized dynamics under a gradient control law to
characterize the system's structural limitations and their dynamic consequences. Our
main contribution is a unified geometric framework for the uncontrollable subspace: an
exact characterization of its rigid-body component and a geometric bound on its
deformational component. We first prove that the uncontrollable rigid-body modes are
exactly the rotations about the actuated node, characterized by the global rotational
subspace $\mathcal{R}_i$. We then introduce the local rotational subspace $\mathcal{T}_i$,
consisting of the motions invisible to the actuator's local measurements, and prove that
for minimally connected actuators, where the actuated node has as many neighbors as the
dimension of the ambient space, the entire uncontrollable subspace is confined to
$\mathcal{T}_i$. Finally, we demonstrate the dynamic implications of this structure by
proving that the ability of the formation to recover its shape is determined by the
alignment of the input with the local component of the rotational rigid-body mode,
directly linking the geometry of hidden modes to disturbance rejection. We illustrate
our results with a case study.
\end{abstract}

\begin{keyword}
Multi-agent systems, Formation control,  Control over networks, Controllability.
\end{keyword}



\end{frontmatter}

\section{Introduction}

The coordination of \ac{mass} is a cornerstone of modern engineering, enabling applications from satellite formation flying to drone swarms and robotic teams (\cite{chung2018survey}). A fundamental task in this domain is formation control, where agents must achieve and maintain a specific geometric pattern. In many real-world scenarios, such as indoor or underwater environments, global positioning is unavailable, forcing agents to rely on local, relative sensing. Distance-based formation control, where the desired shape is defined by a set of inter-agent distances, provides a distributed solution for such settings (\cite{oh2011formation}).

A standard approach to this problem employs a gradient-based control law, where each agent adjusts its motion to minimize the error in its local distance measurements (\cite{krick2009stabilisation}). While the stability and convergence properties of these controllers are well-understood (\cite{SUN201650}), their dynamic performance under external disturbances remains an open area of investigation. How does a localized disturbance, such as a wind gust affecting a single drone, propagate through the formation? How does the placement of the actuator and the geometry of the formation itself affect the system's ability to reject such disturbances? To answer these questions, a formal input-output modeling framework is required.

The analysis of controllability in networked systems has a rich history, particularly for linear consensus-type dynamics. Foundational work has established deep connections between a network's controllability and its underlying graph topology, including the role of graph symmetries (\cite{rahmani2009controllability}) and the placement of leader agents (\cite{tanner2004controllability}). However, these classical results are often tailored to systems with simple, unweighted graph Laplacians. They do not directly address the challenges posed by distance-based formation control, where the system dynamics are governed by a configuration-dependent, weighted Laplacian that introduces geometric constraints and nonlinearities.

This paper develops  a framework to analyze the input-output properties of distance-based formations. We begin by linearizing the nonlinear gradient dynamics around a target equilibrium configuration. This yields a \ac{lti} model whose system matrix is a configuration-dependent, weighted graph Laplacian. This model, while an approximation, allows us to leverage the powerful tools of classical systems theory to analyze the formation's controllability.


The main contributions of this paper are twofold. First, we develop a unified geometric framework for uncontrollable modes. We provide an exact geometric characterization of uncontrollable rigid-body modes via the global rotational subspace $\mathcal{R}_i$, and establish a geometric bound on uncontrollable deformations via the local rotational subspace $\mathcal{T}_i$ for minimally connected actuators. Second, we analyze the dynamic implications of this structure, demonstrating that the system's ability to recover its shape is determined by the alignment of a disturbance with the local component of the standard rotational rigid-body mode. This establishes a direct link between the geometry of hidden modes and disturbance rejection.

The remainder of this paper is organized as follows. Section 2 provides the necessary background on rigidity theory and system-theoretic concepts. Section 3 formally defines the problem, introducing the nonlinear dynamics and the linearized model used for analysis. Section 4 presents our main theoretical contribution: a geometric characterization of the uncontrollable subspace. Finally, Section 5 analyzes the dynamic implications of this geometric structure and illustrates the results with a case study.

\paragraph*{Notations} Throughout, $I_d$ denotes the identity matrix in $\mathbb{R}^d$, {\color{black}$e_\ell$ is the $\ell$-th standard basis vector} in $\mathbb{R}^n$, and $\Omega$ denotes a $d \times d$ skew-symmetric matrix representing an infinitesimal rotation. For $d=2$, {\color{black}every such matrix is a scalar multiple of} $\Omega = \begin{psmallmatrix} 0 & -1 \\ 1 & 0 \end{psmallmatrix}$. The eigenspace associated with an eigenvalue $\lambda$ is denoted by $E_\lambda${\color{black}, and $P_\lambda$ denotes the orthogonal projection onto $E_\lambda$}. The Kronecker product is denoted by $\otimes$. For a stacked vector $v \in \mathbb{R}^{nd}$, its component corresponding to agent $i$ is denoted $v_i \in \mathbb{R}^d$. The degree of a node $i$ in $\Gcal$, i.e., the number of edges incident to $i$, is denoted $\deg(i)$. {\color{black}The set of neighbors of node $i$ is $\mathcal{N}_i := \{j \in \Vcal : (i,j) \in \Ecal\}$, and the centroid of a configuration $p$ is $p_{\mathrm{cm}} := \tfrac{1}{n}\sum_{\ell \in \Vcal} p_\ell$.} The standard Euclidean inner product is denoted by $\langle \cdot, \cdot \rangle$.

\section{Preliminaries}

This section provides the necessary background on graph rigidity and system-theoretic concepts.

\subsection{Rigidity of Graphs}
A formation of $n$ agents is modeled as a \emph{framework} $(\Gcal, p)$, where $\Gcal=(\Vcal, \Ecal)$ is an undirected graph with a vertex set $\Vcal=\{1,\dots,n\}$ and an edge set $\Ecal$ of size $m=|\Ecal|$. The vector $p \in \mathbb{R}^{nd}$ is the stacked vector of agent positions. A reference configuration, used to define the desired inter-agent distances, is denoted $p^*$.

The shape of the framework is encoded by the vector of squared inter-agent distances, known as the \emph{rigidity function}, $r(p) \,:\,\mathbb{R}^{nd} \to \mathbb{R}^m$, with
\begin{equation}\label{rigidity_function}
    r_k(p) = \tfrac{1}{2}\|p_i - p_j\|^2, \quad \text{for each  } k=(i,j) \in \Ecal.
\end{equation}
The first-order change in these distances due to an infinitesimal motion $v \in \mathbb{R}^{nd}$ about the reference configuration $p^*$ is given by the Taylor expansion $r(p^*+tv) - r(p^*) = t R(p^*)v + o(t^2)$, where the Jacobian $R(p^*) := \nabla_p r(p)|_{p=p^*} \in \mathbb{R}^{m\times nd}$ is the \emph{rigidity matrix} evaluated at $p^*$ (\cite{rigidity_of_graphs_1_asimow_roth, rigidity_of_graphs_2_asimow_roth}).


Any distance-preserving first order motion must reside in the nullspace $\ker(R(p^*))$. These allowable motions fall into two distinct physical categories: \ac{rbms}, representing trivial translations and rotations of the entire configuration, and \emph{infinitesimal flexes}, representing internal shape deformations. A framework is infinitesimally rigid if $\ker(R(p^*))$ consists entirely of \ac{rbms}, completely lacking any internal flexes (\cite{connelly2005generic}). For a generic framework in $\mathbb{R}^2$, this is equivalent to the rank condition $\rank(R(p^*))=2n-3$ (\cite{laman1970graphs}).

The rigidity matrix also defines two other fundamental subspaces. The nullspace of its transpose, $\ker(R(p^*)^\top)$, is the space of \emph{self-stresses}. A non-zero self-stress is a set of tensions and compressions on the edges that self-equilibrate at every node. By the fundamental theorem of linear algebra, the image of the transpose, $\image(R(p^*)^\top)$, is the orthogonal complement to the space of RBMs. It therefore represents the space of all infinitesimal \emph{deformations}.


\subsection{Controllability of LTI Systems}
Consider a general \ac{lti} system $\dot{x}=Ax+Bu$. The system's controllability properties determine which internal modes can be influenced by the input.
A key tool for this analysis is the \ac{pbh} test: a mode represented by an eigenvector $v$ of $A$ is \emph{uncontrollable} if and only if it lies in the nullspace of $B^\top$, i.e., $B^\top v = 0$. The set of all such modes spans the uncontrollable subspace $\overline{\mathcal{C}}$, which we refer to as the \emph{hidden} subspace of the system (\cite{Skogestad2005}). 
{\color{black}
For a symmetric state matrix, the \ac{pbh} test yields an explicit modal description of the hidden subspace, which we record for later use.
\begin{lem}[Modal Decomposition of $\overline{\mathcal{C}}$]\label{lem:pbh_decomp}
If $A = A^\top$, then
$\overline{\mathcal{C}} = \bigoplus_{\lambda}\big(E_\lambda \cap \ker(B^\top)\big)$.
\end{lem}
\begin{pf}
A vector $v$ is orthogonal to the controllable subspace if and only if $B^\top A^q v = 0$ for all integers $q \geq 0$. Writing $v = \sum_\lambda P_\lambda v$ gives $B^\top A^q v = \sum_\lambda \lambda^q B^\top P_\lambda v$, and invertibility of the Vandermonde system over the distinct eigenvalues shows that this vanishes for all $q$ if and only if $B^\top P_\lambda v = 0$ for every $\lambda$. \qed
\end{pf}
}

\section{Problem Formulation}

We consider a formation of $n$ agents in $\mathbb{R}^d$ with single-integrator dynamics, $\dot{p}_i = u_i$, a standard model in the distance-based formation control literature (\cite{krick2009stabilisation, oh2011formation, babazadeh2019distance}). The agent interactions are described by an undirected sensing graph $\Gcal$.

The control objective is to steer the formation to a configuration $p$ that satisfies a desired set of inter-agent distances, defined by a reference configuration $p^*$. This is achieved using a gradient-descent control law. We define a potential function based on the collective error in squared edge lengths,
\begin{equation}
    V(p) = \frac{1}{2}\|r(p) - r(p^*)\|^2,
\end{equation}
where $r(p)$ is the rigidity function defined \eqref{rigidity_function}. The control law is the negative gradient of this potential, $\dot{p} = -\nabla_p V(p)$, which is a standard approach in the literature (see, e.g., \cite{krick2009stabilisation, survey_formation_control_2015_automatica}). This yields the closed-loop dynamics:
\begin{equation}
\label{eq:nonlinear_dynamics}
    \dot{p} = -R(p)^\top \big(r(p) - r(p^*)\big).
\end{equation}
Note that $V$ is invariant under continuous rigid-body motions of the formation. Consequently, the equilibrium condition $r(p) = r(p^*)$ does not define an isolated point, but a continuous family of translated and rotated copies of the reference configuration.

To analyze the system's response to a localized exogenous disturbance, we augment these dynamics with a single-channel input applied at one agent $i$,
\begin{equation}
    \dot{p} = -R(p)^\top\big(r(p)-r(p^*)\big) + B w, \label{eq:nonlinear_io_model_ifac}
\end{equation}
where $w \in \mathbb{R}^d$ and $B = e_i \otimes I_d$. The single-actuator setting is the simplest abstraction in which to expose how actuator placement and graph geometry interact, and serves as a building block for multi-input extensions.

Direct analysis of the nonlinear model \eqref{eq:nonlinear_io_model_ifac} is challenging. To gain tractable insights, we linearize the system around the equilibrium point $(p^*, w=0)$. Letting {\color{black}$x := p - p^*$ denote the state deviation}, the first-order Taylor expansion of the dynamics yields the \ac{lti} model:
\begin{equation}
    \;{\color{black}\dot{x} = A\,x + B\,w},\;
    \label{eq:first-order-lti-ifac}
\end{equation}
where the system matrix $A$ is the weighted graph Laplacian, or ``stiffness matrix," given by:
\begin{equation}
    A = -R(p^*)^\top R(p^*).
\end{equation}
The analysis that follows characterizes the input-output behavior of \eqref{eq:first-order-lti-ifac} in a neighborhood of $p^*$ where this linearization is valid; all subsequent statements, including the shape-recovery dichotomy of Section~5, are properties of this linearized model. Future work will focus on expanding this geometric interpretation to the full nonlinear dynamics.

\section{Geometric Characterization\\ of Hidden Modes}


This section develops a geometric characterization of the hidden modes for the linearized system \eqref{eq:first-order-lti-ifac}. We first establish their existence and algebraic structure before revealing their geometric form.

\subsection{Algebraic Foundations of Hidden Modes}

We begin by establishing that for an infinitesimally rigid formation, hidden modes are an unavoidable consequence of using a single actuator.


\begin{prop}[Existence of Hidden RBMs]
\label{prop:existence_hidden_rbm}
Consider the linearized system \eqref{eq:first-order-lti-ifac} for an infinitesimally rigid framework in
$\mathbb{R}^d$ with $d\ge2$, actuated at a single node $i$. Then
\[
  \dim(\overline{\mathcal{C}}\cap E_0) = \frac{d(d-1)}{2} \;\ge\; 1.
\]
\end{prop}
\begin{pf}
For an infinitesimally rigid framework, $\dim(E_0)=d(d+1)/2$. By
Lemma~\ref{lem:pbh_decomp}, $\overline{\mathcal{C}}\cap E_0 = E_0\cap\ker(B^\top)$.
Since $B=e_i\otimes I_d$, $\dim(\ker B^\top)=dn-d$, and $E_0+\ker B^\top=\mathbb{R}^{dn}$:
any $w\in\operatorname{Im}B\cap E_0^\perp$ is supported only at node $i$ and orthogonal
to the $d$ translations, forcing $w=0$. Grassmann's identity then gives
$\dim(E_0\cap\ker B^\top) = \frac{d(d+1)}{2}+(dn-d)-dn = \frac{d(d-1)}{2}$, which is at
least $1$ for $d\ge2$. \qed
\end{pf}

Having established their existence, we now seek to understand their structure.

{\color{black}
\begin{prop}[Uncontrollable Subspace Decomposition]\label{prop:decomp}
The uncontrollable subspace splits into its rigid-body and deformational
components as
\begin{equation}\label{eq:decomp}
  \bar{\mathcal{C}}
  = \underbrace{\big(\ker R(p^*)\cap\ker B^\top\big)}_{\text{hidden RBMs}}
  \;\oplus\;
  \underbrace{\bigoplus_{\lambda\neq0}\big(E_\lambda\cap\ker B^\top\big)}_{\text{hidden deformations}} ,
\end{equation}
and the deformational term admits the equivalent form
$\bigoplus_{\lambda\neq0}(E_\lambda\cap\ker B^\top)
 = \bar{\mathcal{C}}\cap\image R(p^*)^\top$.
\end{prop}

\begin{pf}
{\color{black}By Lemma~\ref{lem:pbh_decomp}}, $\bar{\mathcal{C}}=\bigoplus_\lambda(E_\lambda\cap\ker B^\top)$. Isolating $\lambda=0$ and using
$E_0=\ker R(p^*)$ yields \eqref{eq:decomp}.

For the equivalent form, since $A$ is symmetric, $\mathbb{R}^{dn}=\bigoplus_\lambda E_\lambda$ and every $v$ admits
the unique decomposition $v=\sum_\lambda v_\lambda$ with $v_\lambda=P_\lambda v\in E_\lambda$.

(i) If $v\in\bar{\mathcal{C}}$, then $v_\lambda\in\ker B^\top$ for every $\lambda$: indeed
$v=\sum_\lambda w_\lambda$ for some $w_\lambda\in E_\lambda\cap\ker B^\top$, and uniqueness
forces $w_\lambda=v_\lambda$.

(ii) Since $\image R(p^*)^\top=(\ker R(p^*))^\perp=E_0^{\perp}$, we have
$v\in\image R(p^*)^\top\iff P_0v=0\iff v_0=0$.

Combining (i) and (ii), $v\in\bar{\mathcal{C}}\cap\image R(p^*)^\top$ if and only if
$v=\sum_{\lambda\neq0}v_\lambda$ with $v_\lambda\in E_\lambda\cap\ker B^\top$. \qed
\end{pf}
}

This decomposition provides the algebraic foundation. To find a geometric interpretation for each term in \eqref{eq:decomp}, we apply the \ac{pbh} test to our single-actuator system. With $B = e_i \otimes I_d$, the uncontrollability condition $B^\top v = 0$ simplifies to a direct ``pinning'' constraint at the agent level.


\begin{lem}[Pinning Constraint]\label{lem:pinning}
An eigenvector $v$ is uncontrollable from an input applied at agent $i$ if and only if its component at that agent vanishes, i.e., $v_i = 0$.
\end{lem}
\begin{pf}
From the \ac{pbh} test, $v$ is uncontrollable iff $B^\top v = (e_i^\top \otimes I_d) v = v_i = 0$. \qed
\end{pf}

\subsection{Geometric Characterization of Hidden RBMs}
\label{subsec:hidden_rbms}

We now apply the pinning constraint of Lemma \ref{lem:pinning} to give the two components of the decomposition \eqref{eq:decomp} their geometric form, beginning with the rigid-body part: what does the hidden \ac{rbm} guaranteed by Proposition \ref{prop:existence_hidden_rbm} look like?



\begin{prop}[Uncontrollable RBM is a Pure Rotation] \label{prop:R_i}
For some skew-symmetric $\Omega \in \mathbb{R}^{d\times d}$, the uncontrollable \ac{rbm}
subspace $\ker R(p^*)\cap\ker B^\top$ equals the $d(d-1)/2$-dimensional rotational
subspace
\[
  \mathcal{R}_i(p^*) := \big\{\, v \in \mathbb{R}^{dn} \mid v_\ell = \Omega(p_\ell^* - p_i^*),\ \forall \ell \in \mathcal{V} \,\big\}.
\]
\end{prop}

\begin{pf}
Any \ac{rbm} admits a unique decomposition {\color{black}$v_\ell = \xi + \Omega(p_\ell^*-p_{\mathrm{cm}})$ into a translation $\xi \in \mathbb{R}^d$} and a rotation about the centroid. The pinning constraint $v_i=0$ from Lemma~\ref{lem:pinning} fixes {\color{black}$\xi=-\Omega(p_i^*-p_{\mathrm{cm}})$}, which after substitution yields {\color{black}$v_\ell=\Omega(p_\ell^*-p_i^*)$}, a pure rotation about node $i$. \qed
\end{pf}



This geometric interpretation is illustrated in Figure~\ref{fig:global_rotation_illustration}. This proof reveals a key geometric insight: the algebraic ``pinning constraint" $v_i=0$ forces the center of rotation of any hidden \ac{rbm} to coincide with the actuated node $i$. This establishes a direct and intuitive link between the placement of the actuator and the geometric form of the resulting uncontrollable mode. 


\begin{figure*}[t]
    \begin{minipage}[t]{0.47\linewidth}
        \centering
        \begin{tikzpicture}[scale=0.95,
  vtx/.style={circle, draw=black, fill=white, inner sep=2.0pt, font=\scriptsize},
  focus/.style={circle, draw=black, fill=red!70, inner sep=2.2pt, font=\scriptsize\bfseries, line width=0.8pt},
  starEdge/.style={draw=black, line width=1.3pt},
  otherEdge/.style={draw=gray!55, dashed, line width=1.0pt},
  armStyle/.style={draw=red!60, dash pattern=on 1.2mm off 0.8mm, line width=0.9pt},
  tanStyle/.style={draw=red!80!black, line width=1.4pt, -{Latex[length=1.6mm]}, line cap=round}
]
  \useasboundingbox (-2.8,-2.6) rectangle (2.9,2.4);
  \coordinate (ni) at (0,0);
  \coordinate (na) at (2.0, 0.9);
  \coordinate (nb) at (1.2,-1.1);
  \coordinate (nc) at (-1.9, 0.5);
  \coordinate (nd) at (-0.8,-1.7);
  \draw[otherEdge] (na)--(nb);
  \draw[otherEdge] (na)--(nd);
  \draw[otherEdge] (nb)--(nd);
  \draw[otherEdge] (nc)--(nd);
  \draw[starEdge] (ni)--(na);
  \draw[starEdge] (ni)--(nb);
  \draw[starEdge] (ni)--(nc);
  \foreach \v in {na,nb,nc,nd}{
    \draw[armStyle] (ni)--(\v);
    \draw[tanStyle]
      ($(\v)!1.5mm!90:(ni)$) --
      ($($(\v)!1.5mm!90:(ni)$)!4.5mm!($(\v)!2.5mm!90:(ni)$)$);
  }
  \node[focus] at (ni) {$i$};
  \node[vtx]   at (na) {$\alpha$};
  \node[vtx]   at (nb) {$\beta$};
  \node[vtx]   at (nc) {$\gamma$};
  \node[vtx]   at (nd) {$\ell$};
\end{tikzpicture}
        \caption{Geometric interpretation of the uncontrollable rotational subspace $\mathcal{R}_i(p^*)$. The red arrows depict a pure rotation of the \emph{entire} framework about the actuated node $i$ (red filled). Dashed red lines indicate the rotation arms from $i$ to every other node, showing that this global rigid-body motion is generated jointly by all relative positions with respect to $i$. The dashed gray edges show the underlying sensing graph and the solid black edges form the incident edges.}
        \label{fig:global_rotation_illustration}
    \end{minipage}
    \hfill
    \begin{minipage}[t]{0.47\linewidth}
        \centering
        \begin{tikzpicture}[scale=0.95,
  vtx/.style={circle, draw=black, fill=white, inner sep=2.0pt, font=\scriptsize},
  focus/.style={circle, draw=black, fill=red!70, inner sep=2.2pt, font=\scriptsize\bfseries, line width=0.8pt},
  starEdge/.style={draw=black, line width=1.3pt},
  otherEdge/.style={draw=gray!55, dashed, line width=1.0pt},
  tanStyle/.style={draw=red!80!black, line width=1.4pt, -{Latex[length=1.6mm]}, line cap=round}
]
  \useasboundingbox (-2.8,-2.6) rectangle (2.9,2.4);
  \coordinate (ni) at (0,0);
  \coordinate (na) at (2.0, 0.9);
  \coordinate (nb) at (1.2,-1.1);
  \coordinate (nc) at (-1.9, 0.5);
  \coordinate (nd) at (-0.8,-1.7);

  \draw[otherEdge] (na)--(nb);
  \draw[otherEdge] (na)--(nd);
  \draw[otherEdge] (nb)--(nd);
  \draw[otherEdge] (nc)--(nd);

  \draw[starEdge] (ni)--(na);
  \draw[starEdge] (ni)--(nb);
  \draw[starEdge] (ni)--(nc);

  \foreach \v/\side/\len in {na/90/5.5mm, nb/-90/4.0mm, nc/90/6.5mm}{
    \draw[tanStyle]
      ($(\v)!1.5mm!\side:(ni)$) --
      ($($(\v)!1.5mm!\side:(ni)$)!\len!($(\v)!2.5mm!\side:(ni)$)$);
  }

  \node[focus] at (ni) {$i$};
  \node[vtx]   at (na) {$\alpha$};
  \node[vtx]   at (nb) {$\beta$};
  \node[vtx]   at (nc) {$\gamma$};
  \node[vtx]   at (nd) {$\ell$};
\end{tikzpicture}
        \caption{The local rotational subspace $\mathcal{T}_i$ at the actuated node $i$ (red). Solid black edges form the incident star; dashed gray edges are non-incident. The red arrows depict elementary motions $\tau_{i\to j}$, where each neighbor {\color{black}$j\in\{\alpha,\beta,\gamma\}$} can rotate independently about $i$. Because these local rotations are decoupled, they generate internal shape deformations along the non-incident edges. These local motions only manifest as true dynamic hidden modes if they simultaneously satisfy the global eigenvector conditions of $\overline{\mathcal{C}}$.}
        \label{fig:local_rotation_illustration}
    \end{minipage}
\end{figure*}

\subsection{Locally Hidden Deformations}
\label{subsec:hidden_deformations}

Having characterized the hidden \ac{rbm}s, we now generalize our approach to find hidden deformations. The gradient controller at node $i$ is driven only by errors in its incident edges. A motion is therefore ``locally hidden" from the controller at node $i$ if it generates no first-order change in these local edge lengths.

To formalize this, we consider the components of the rigidity function, $r_k(p) = \tfrac{1}{2}\|p_i - p_j\|^2$, corresponding to the edges incident to node $i$. These form the \emph{star relations} at $i$. The first-order change in these relations due to an infinitesimal motion $v$ is given by
\[
\delta r_{ij} = (p_i^*-p_j^*)^\top(v_i-v_j), \; \text{for each edge } (i,j) \in \Ecal.
\]
A motion $v$ is locally hidden from node $i$ if it produces no change in these local measurements, i.e., {\color{black}$\delta r_{ij}=0$ for all neighbors $j \in \mathcal{N}_i$}. This motivates the following definition, which adds the pinning constraint $v_i=0$.

\begin{defn}[The Local Rotational Subspace $\mathcal{T}_i(\mathcal{G},p^*)$] \label{def:T_i}
The set of infinitesimal motions that fix node $i$ and preserve the lengths of all its
incident edges to first order is
\[
  \mathcal{T}_i := \{v \in \mathbb{R}^{dn} : v_i = 0 \text{ and } (p_j^* - p_i^*)^\top v_j = 0,\ \forall (i,j) \in \mathcal{E} \}.
\]
\end{defn}

To compute the dimension of $\mathcal{T}_i$, we construct an explicit basis by counting the constraints in Definition~\ref{def:T_i}. The pinning condition $v_i=0$ removes $d$ degrees of freedom from the ambient $\mathbb{R}^{dn}$. The remaining $d(n-1)$ dimensions split according to the actuator's neighborhood:
\begin{enumerate}
    \item [(i)] For each neighbor $j$ of $i$, the length-preserving condition $(p_j^*-p_i^*)^\top v_j = 0$ removes one additional degree of freedom, confining $v_j$ to the $(d-1)$-dimensional hyperplane orthogonal to the incident edge. In the planar case ($d=2$), this hyperplane reduces to the line spanned by the elementary rotation
    \begin{equation}
        \tau_{i\to j} := (e_j \otimes I_d)\,\Omega\,(p_j^*-p_i^*),
    \end{equation}
    which rotates only neighbor $j$ infinitesimally about $i$.
    \item [(ii)] Each non-neighbor $\ell$ is unconstrained by Definition~\ref{def:T_i} and contributes its full $d$ standard basis directions.
\end{enumerate}
Summing these contributions yields $$\dim(\mathcal{T}_i) = d(n-1) - \deg(i).$$ 
By construction, every motion in $\mathcal{T}_i$ preserves the edges incident to $i$ and is therefore locally hidden to the gradient law at $i$. At the same time, such motions are not globally rigid: a neighbor $j$ rotating about $i$ leaves the edge $(i,j)$ intact but generally changes all distances from $j$ to other nodes, producing a shape-changing deformation that goes undetected by the actuator. These elementary motions are illustrated in Figure~\ref{fig:local_rotation_illustration}.

\begin{prop}[Geometric Subspace Inclusion] \label{prop:inclusion_R_i_T_i}
The rotational subspace $\mathcal{R}_i(p^*)$ is a subspace of the local rotational subspace $\mathcal{T}_i(\Gcal, p^*)$, i.e., $\mathcal{R}_i(p^*) \subseteq \mathcal{T}_i(\Gcal, p^*)$.
\end{prop}
\begin{pf}
For $v\in\mathcal{R}_i(p^*)$, $v_i=\Omega(p_i^*-p_i^*)=0$ satisfies the pinning constraint, and $(p_j^*-p_i^*)^\top \Omega(p_j^*-p_i^*)=0$ by skew-symmetry of $\Omega$. \qed
\end{pf}

The inclusion highlights a key distinction: $\mathcal{R}_i$ contains only the single global RBM that is a pure rotation about $i$, whereas $\mathcal{T}_i$ also contains the elementary motions $\tau_{i\to j}$, which are locally rigid at $i$ but globally deformational and lie in $\image(R(p^*)^\top)$.

\begin{thm}[Geometric Bound for $\deg(i)=d$] \label{thm:geometric_bound}
For an infinitesimally rigid framework $(\mathcal{G}, p^*)$ actuated at a node $i$ with minimal degree $\deg(i) = d$, {\color{black}whose incident edge vectors $\{p_j^* - p_i^*\}_{j \in \mathcal{N}_i}$ are linearly independent (as holds for generic configurations),} the uncontrollable subspace is geometrically bounded by the local rotational subspace, i.e.,
$$ \overline{\mathcal{C}} \subseteq \mathcal{T}_i(\mathcal{G}, p^*). $$
\end{thm}
\begin{pf}
{\color{black}Let $v \in \overline{\mathcal{C}}$. By part (i) in the proof of Proposition~\ref{prop:decomp}, $v = \sum_\lambda v_\lambda$, where each $v_\lambda \in E_\lambda \cap \ker B^\top$ is an uncontrollable eigenvector. Since $\mathcal{T}_i$ is a subspace, it suffices to show that each eigencomponent lies in $\mathcal{T}_i$; we may therefore assume that $v$ itself is an uncontrollable eigenvector, satisfying the pinning constraint $v_i = 0$ (Lemma~\ref{lem:pinning}) and the eigenvector equation $-R(p^*)^\top R(p^*)v = \lambda v$.}

Evaluating the $i$-th vector block of the eigenvector equation yields the equilibrium condition at the actuator:
\begin{equation} \label{eq:equilibrium}
\sum_{j \in \mathcal{N}_i} (p_i^* - p_j^*) \delta r_{ij} = 0,
\end{equation}
where $\delta r_{ij} = -(p_i^* - p_j^*)^\top v_j$ is the first-order change in the
length of the incident edge $(i,j)$ given $v_i=0$.


Since the framework is infinitesimally rigid, the rank condition of the rigidity
matrix inherently requires $\deg(i) \ge d$ to prevent local infinitesimal flexes.
By our premise, $\deg(i) = d$. {\color{black}Since the $d$ incident edge vectors are linearly
independent by assumption}, they form a basis for
$\mathbb{R}^d$, and the equilibrium condition \eqref{eq:equilibrium} is satisfied if
and only if $\delta r_{ij} = 0$ for all $j \in \mathcal{N}_i$.

Consequently, the mode $v$ satisfies both the pinning constraint ($v_i=0$) and the preservation of all incident edges ($\delta r_{ij} = 0$). These are exactly the defining geometric conditions of Definition~\ref{def:T_i}, implying $v \in \mathcal{T}_i$. Therefore, the geometric inclusion $\overline{\mathcal{C}} \subseteq \mathcal{T}_i$ holds. \qed
\end{pf}


By substituting the geometric bound from Theorem \ref{thm:geometric_bound} into the algebraic decomposition of Proposition \ref{prop:decomp}, and accounting for the rigid-body inclusion ($\mathcal{R}_i \subset \mathcal{T}_i$) from Proposition \ref{prop:inclusion_R_i_T_i}, we immediately isolate the deformational hidden modes.

\begin{cor}\label{cor:bound}
{\color{black}Under the assumptions of Theorem~\ref{thm:geometric_bound},}
\begin{equation}\label{eq:refined_bound}
  \bar{\mathcal{C}} \subseteq
  \mathcal{R}_i(p^*)\oplus\big(\mathcal{T}_i(\mathcal{G},p^*)\cap\image R(p^*)^\top\big)
  \subseteq \mathcal{T}_i(\mathcal{G},p^*).
\end{equation}
\end{cor}

\begin{pf}
By Proposition~\ref{prop:decomp}, $\bar{\mathcal{C}}=(\bar{\mathcal{C}}\cap\ker R(p^*))
\oplus(\bar{\mathcal{C}}\cap\image R(p^*)^\top)$. The first term is $\ker R(p^*)\cap
\ker B^\top=\mathcal{R}_i(p^*)$ by Proposition~\ref{prop:R_i}. For the second,
Theorem~\ref{thm:geometric_bound} gives $\bar{\mathcal{C}}\subseteq\mathcal{T}_i$, hence
$\bar{\mathcal{C}}\cap\image R(p^*)^\top\subseteq\mathcal{T}_i\cap\image R(p^*)^\top$.
The final inclusion holds since both summands lie in $\mathcal{T}_i$, the first by
Proposition~\ref{prop:inclusion_R_i_T_i}. \qed
\end{pf}


{\color{black}
\begin{rem}[Dimension of the bound~\eqref{eq:refined_bound}]
Let $\deg(i)=d$ and $n\ge d+2$, so the graph contains a node $\ell$ that is neither $i$
nor a neighbor of $i$. Every motion in $\image(e_\ell\otimes I_d)$ then satisfies
Definition~\ref{def:T_i}, since it vanishes at $i$ and at every neighbor of
$i$; such motions are not rigid, yet they are not orthogonal to the translations.
Together with $\mathcal{R}_i\subseteq\mathcal{T}_i$
(Proposition~\ref{prop:inclusion_R_i_T_i}), the orthogonal projection $P_0$ onto
$\ker R(p^*)$ therefore satisfies $P_0\mathcal{T}_i=\ker R(p^*)$, a rank condition,
which holds for generic $p^*$. {\color{black}Indeed, since $\ker R(p^*)$ is spanned by
$\mathcal{R}_i$ and $d$ translations, this condition only concerns
the translational components of the motions above.} Equivalently $\mathcal{T}_i+\image R(p^*)^\top
=\mathbb{R}^{dn}$, and Grassmann's formula gives
\begin{align*}
    \dim\big(\mathcal{T}_i\cap\image R(p^*)^\top\big)
  &=\dim(\mathcal{T}_i)+\dim\big(\image R(p^*)^\top\big)-dn\\
  &=d(n-2)-\tfrac{d(d+1)}{2}.
\end{align*}
With $\dim(\mathcal{R}_i)=d(d-1)/2$, Corollary~\ref{cor:bound} then yields
\[
  \tfrac{d(d-1)}{2}\;\le\;\dim(\bar{\mathcal{C}})\;\le\;d(n-3),
\]
where the lower bound holds unconditionally, since
$\dim(\overline{\mathcal{C}}) \ge \dim(\overline{\mathcal{C}}\cap E_0) = \tfrac{d(d-1)}{2}$
by Proposition~\ref{prop:existence_hidden_rbm}, whereas the upper bound requires
$\deg(i)=d$, $n\ge d+2$, and the rank condition above. The bounding
subspace of Corollary~\ref{cor:bound} is thus a proper subspace of $\mathcal{T}_i$, of
dimension $d$ less than $\dim(\mathcal{T}_i)=d(n-2)$.
\end{rem}
}

\section{Dynamic Implications of Hidden Modes}

We now analyze the dynamic consequences of the hidden modes, investigating how a localized input at agent $i$ affects the formation's ability to recover its shape. This reveals a direct link between the geometry of the uncontrollable subspace $\mathcal{R}_i$ and the system's steady-state response. For clarity of exposition, we take $d=2$ throughout this section; the general case follows by replacing the rotational mode $v_r$ below with a basis of rotational modes.

Consider an impulsive input $w(t) = w_0 \delta(t)$, where $w_0 \in \mathbb{R}^2$ is the direction of the impulse. The final state deviation is given by the projection of the initial impulse onto the
\ac{rbm} eigenspace $E_0 = \ker(R(p^*))$. The impulse sets the initial
condition $x(0^+) = Bw_0$; since $A$ is \ac{nsd} with $\ker(A) = E_0$, the components of
$x(0^+)$ along the strictly negative eigenspaces decay as $t\to\infty$, leaving
\begin{equation}
    \lim_{t\to\infty} x(t) = P_0 B w_0.
\end{equation}
Let $\{v_x, v_y, v_r\}$ be an orthonormal basis for the \ac{rbm} space, representing x-translation, y-translation, and rotation about the center of mass, respectively. The final state will be a linear combination of these modes,
\begin{equation}
\label{eq:steady_state_combination}
    \lim_{t\to\infty} x(t) = c_x v_x + c_y v_y + c_r v_r,
\end{equation}
where the coefficients are determined by the input: 
\[c_x = \langle v_x, Bw_0 \rangle, \; c_y = \langle v_y, Bw_0 \rangle, \; c_r = \langle v_r, Bw_0 \rangle.\]

The key insight is that shape recovery depends entirely on whether the global rotational
mode $v_r$ is excited. Since $B = e_i\otimes I_d$, the excitation reduces to a purely
local condition: $c_r = \langle v_r, Bw_0\rangle = \langle [v_r]_i, w_0\rangle$, where
$[v_r]_i \propto \Omega(p_i^* - p_{cm})$ is the velocity of the centroidal rotation at
the actuated node.

\begin{prop}[Shape Recovery]
    \label{prop:shape_recovery_dichotomy}
Consider an impulsive input $w(t) = w_0 \delta(t)$ applied at agent $i$. The ability of the formation to recover its shape is determined by the alignment between the input vector $w_0$ and the local rotational vector $[v_r]_i$.
\begin{enumerate}
    \item[(i)] \emph{Perfect Shape Recovery:} \\If the input is orthogonal to the local rotational direction, $\langle [v_r]_i, w_0 \rangle = 0$, then the rotational mode is not excited ($c_r=0$). The resulting steady-state motion is a pure translation of the entire framework. All inter-agent distances are preserved, and the steady-state edge error is zero.
    
    \item[(ii)] \emph{Persistent Shape Distortion:} \\If the input has a component along the local rotational direction, $\langle [v_r]_i, w_0 \rangle \neq 0$, then the rotational mode is excited ($c_r \neq 0$). The resulting steady-state motion is a combination of translation and rotation, which causes a persistent change in the inter-agent distances and a non-zero steady-state edge error.
\end{enumerate}
\end{prop}


\begin{pf}
The coefficient of the rotational mode is
$c_r = \langle v_r, Bw_0 \rangle = \langle v_r, (e_i \otimes I_d)w_0 \rangle
     = \langle [v_r]_i, w_0 \rangle$.
We analyze the steady-state edge error at the displaced configuration
$p' = p^* + x(\infty)$: an edge $(j,\ell)$ is preserved if
$\|p'_j - p'_\ell\| = \|p_j^* - p_\ell^*\|$.

\emph{Case 1 (Perfect shape recovery).}
If $c_r = 0$, then $x(\infty) = c_x v_x + c_y v_y$ is a pure translation: there exists
$\xi \in \mathbb{R}^d$ such that $[x(\infty)]_\ell = \xi$ for all $\ell \in \Vcal$. For
every edge $(j,\ell)$,
\[
  p'_j - p'_\ell = (p_j^* + \xi) - (p_\ell^* + \xi) = p_j^* - p_\ell^*,
\]
so all relative positions, and hence all edge lengths, are preserved, and the
steady-state edge error is zero.

\emph{Case 2 (Persistent shape distortion).}
If $c_r \neq 0$, the steady state contains a rotational component:
$[x(\infty)]_\ell = \xi + c_r\,\Omega(p_\ell^* - p_{cm})$ for all $\ell \in \Vcal$.
The translation $\xi$ and the terms in $p_{cm}$ cancel in relative positions, so for
any edge $(j,\ell)$, with edge vector $z := p_j^* - p_\ell^*$,
\begin{align*}
  p'_j - p'_\ell
  &= z + c_r\,\Omega\big(p_j^* - p_\ell^*\big)
   = (I_d + c_r\Omega)\,z.
\end{align*}
Using $\Omega^\top = -\Omega$, the new squared length is
\begin{align*}
  \|(I_d + c_r\Omega)z\|^2
  &= z^\top (I_d - c_r\Omega)(I_d + c_r\Omega)\, z \\
  &= z^\top (I_d - c_r^2\,\Omega^2)\, z
   = \|z\|^2 + c_r^2\,\|\Omega z\|^2 .
\end{align*}
For $d=2$, any nonzero skew-symmetric $\Omega$ is invertible, so
$\|\Omega z\|^2 > 0$ for every nonzero edge vector $z$. Hence for any $c_r \neq 0$
the squared length strictly increases on every edge, and the edge error is persistent.
\qed
\end{pf}

This result reveals a fundamental dichotomy in the system's response, illustrated in
Figures~\ref{fig:shape_recovery_dichotomy_1} and~\ref{fig:shape_recovery_dichotomy_2}
and discussed in the following subsection.

\subsection{Case Study: Minimally Rigid Framework}

We consider a minimally rigid formation of four agents in the plane, with an impulsive
input applied at node~1. The simulations integrate the linearized model
\eqref{eq:first-order-lti-ifac}.

\begin{figure}[!ht]
    \centering
    \begin{subfigure}{0.24\textwidth}
        \centering
        \includegraphics[width=\linewidth]{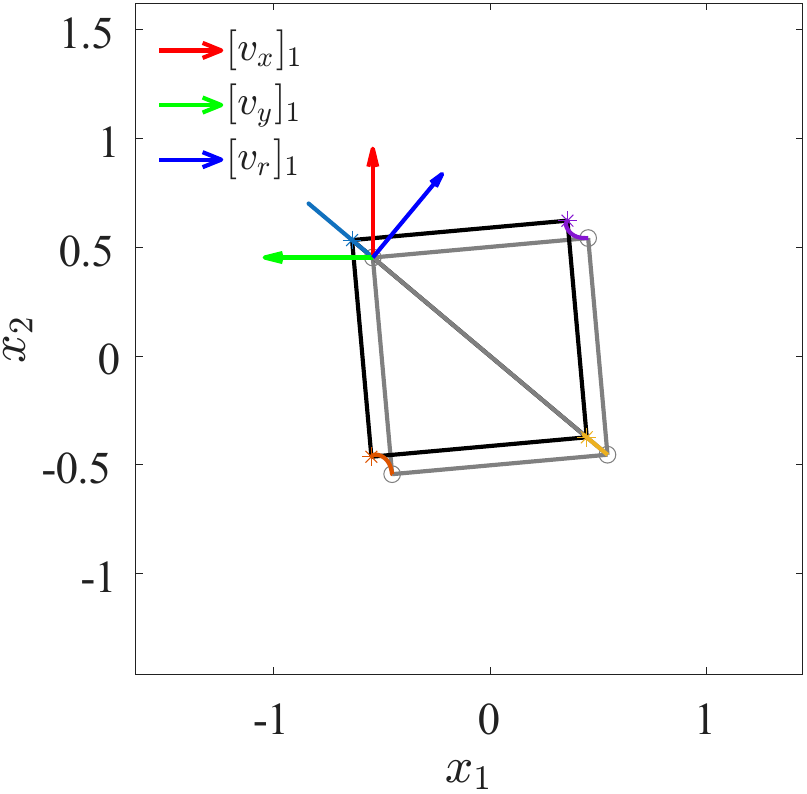}
        \caption{}
        \label{fig:act_dir_1}
    \end{subfigure}
    \begin{subfigure}{0.24\textwidth}
        \centering
        \includegraphics[width=1\linewidth]{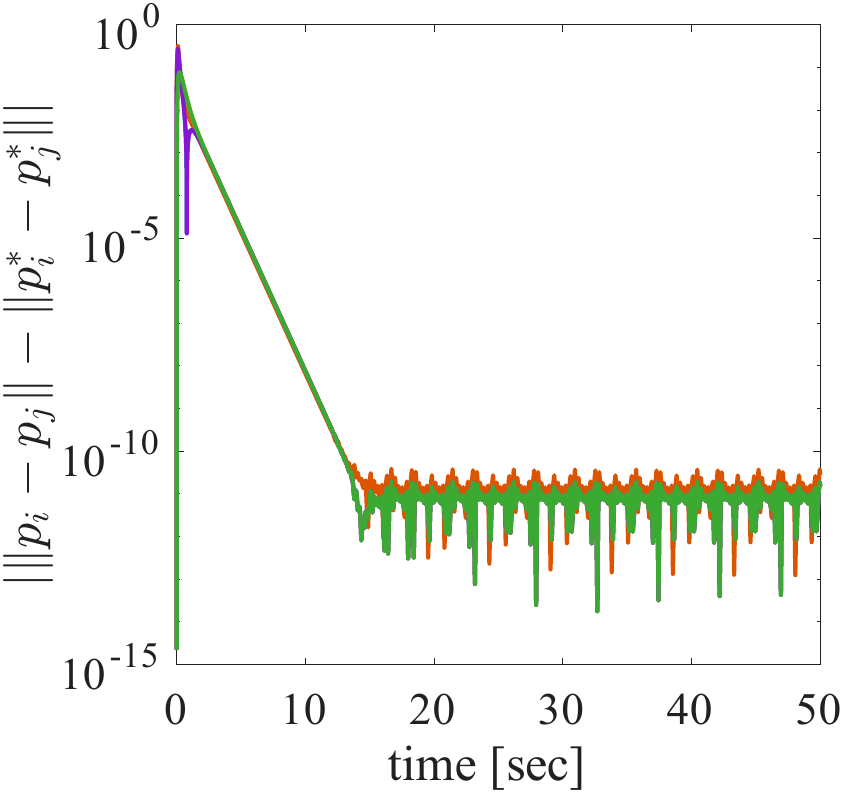}
        \caption{}
        \label{fig:act_dir_1_edge_err}
    \end{subfigure}
    \caption{Shape recovery with an orthogonal input. (a) An impulse orthogonal to the local rotational RBM direction at the actuated node results in a pure translation. (b) The edge length errors decay to zero.}
    \label{fig:shape_recovery_dichotomy_1}
\end{figure}

First, we consider an input $w_0$ orthogonal to the local rotational vector at the actuated node, $\langle [v_r]_1, w_0 \rangle = 0$. This corresponds to Case~1 of Proposition~\ref{prop:shape_recovery_dichotomy}, which predicts that the rotational mode is not excited ($c_r=0$). Figure~\ref{fig:shape_recovery_dichotomy_1}(a) shows the resulting motion: the formation undergoes a pure translation. Figure~\ref{fig:shape_recovery_dichotomy_1}(b) confirms this, showing that all edge length errors decay to zero after the initial transient, and the formation perfectly recovers its shape.

\begin{figure}[!h]
    \centering
    \begin{subfigure}{0.24\textwidth}
        \centering
        \includegraphics[width=\linewidth]{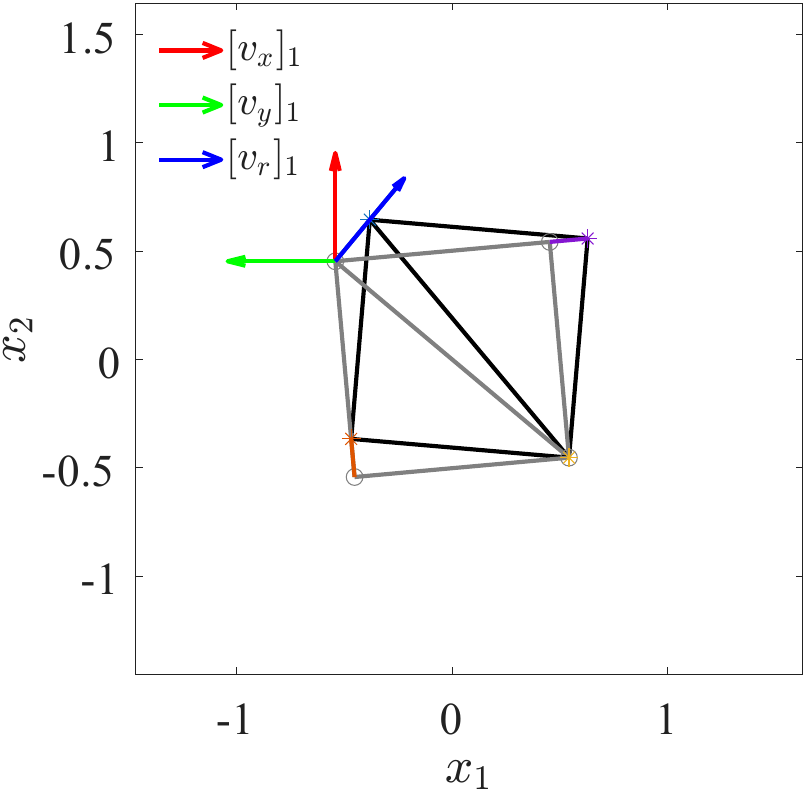}
        \caption{}
        \label{fig:act_dir_2}
    \end{subfigure}
    \begin{subfigure}{0.24\textwidth}
        \centering
        \includegraphics[width=1.07\linewidth]{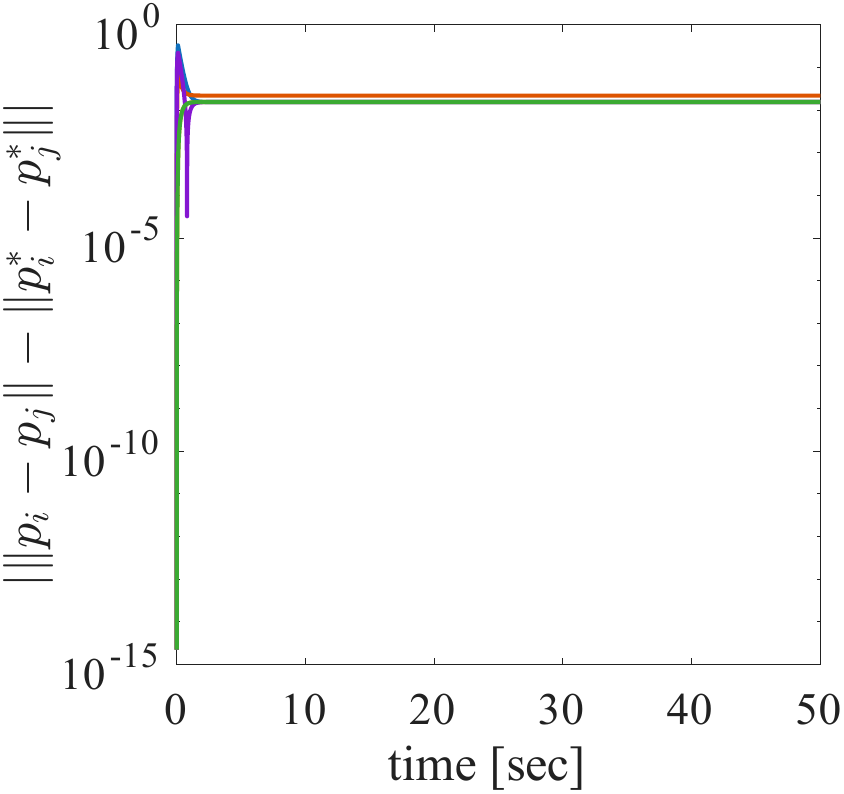}
        \caption{}
        \label{fig:act_dir_2_edge_err}
    \end{subfigure}
    \caption{Shape distortion with an aligned input. (a) An impulse with a component along the local rotational RBM direction excites a rotational motion. (b) The formation deforms, resulting in persistent edge length errors.}
    \label{fig:shape_recovery_dichotomy_2}
\end{figure}

Next, we consider an input $w_0$ with a component along the local rotational vector, $\langle [v_r]_1, w_0 \rangle \neq 0$ (Case~2). This input excites the rotational mode ($c_r\neq 0$): Figure~\ref{fig:shape_recovery_dichotomy_2}(a) shows that the final motion combines a
translation with a rotation of the whole formation. Unlike a translation, this
displacement does not preserve the inter-agent distances, and
Figure~\ref{fig:shape_recovery_dichotomy_2}(b) shows the edge length errors converge to
persistent, non-zero values, providing a clear dynamic validation of our geometric
analysis.


The dichotomy admits a clean geometric interpretation in \ac{rbm} coordinate space.
The steady-state map $w_0 \mapsto (c_x, c_y, c_r)$, with $c_x = \langle [v_x]_i, w_0
\rangle$, $c_y = \langle [v_y]_i, w_0 \rangle$, $c_r = \langle [v_r]_i, w_0 \rangle$,
is linear in $w_0 \in \mathbb{R}^2$, so the reachable \ac{rbm} coordinates form the
two-dimensional ``controllable plane'' shown in
Figure~\ref{fig:unreachable_rbm_illustration}. Its normal,
$n_c=\begin{psmallmatrix}-[v_r]_i\\1\end{psmallmatrix}$, is precisely the coordinate
representation of the hidden rotation $\mathcal{R}_i(p^*)$: substituting these
coordinates into \eqref{eq:steady_state_combination} recovers the pure rotation about
the actuated node. Inputs orthogonal to $[v_r]_i$ further restrict $c_r=0$, placing the
outcome in the translation subspace and yielding perfect shape recovery.

\begin{figure}[!ht]
    \centering
    \tdplotsetmaincoords{85}{120}

\begin{tikzpicture}[
    tdplot_main_coords, 
    scale=2.2, 
    font=\sffamily
]

    \draw[-{Stealth[length=2.5mm]}] (0,0,0) -- (-1.8,0,0) node[anchor=north]{$c_x$};
    \draw[-{Stealth[length=2.5mm]}] (0,0,0) -- (0,1.5,0) node[anchor=west]{$c_y$};
    \draw[-{Stealth[length=2.5mm]}] (0,0,0) -- (0,0,1.2) node[anchor=south]{$c_r$};

    \draw[fill=gray!15, fill opacity=0.3, draw=gray!40] 
        (-1.4,-1.1,0) -- (1.4,-1.1,0) -- (1.4,1.1,0) -- (-1.4,1.1,0) -- cycle;
    
    \def\nx{0.6} \def\ny{-0.3} \def\nz{0.8}
    \coordinate (n) at (\nx, \ny, \nz);
    
    \coordinate (b1) at (\ny, -\nx, 0); 
    \coordinate (b2) at ({\nx*\nz},{\ny*\nz},{-(\nx*\nx + \ny*\ny)}); 
    
    \draw[fill=blue, fill opacity=0.3] 
        ($-1.2*(b1)$) -- 
        ($1.2*(b1)$) -- 
        ($1.2*(b1) + 0.5*(b2)$) -- 
        ($-1.2*(b1) + 0.5*(b2)$) -- cycle;
    \draw[fill=blue, fill opacity=0.6] 
        ($-1.2*(b1) - 0.5*(b2)$) -- 
        ($1.2*(b1) - 0.5*(b2)$) -- 
        ($1.2*(b1)$) -- 
        ($-1.2*(b1)$) -- cycle;

    \draw[orange!85!black, very thick] ($-1.2*(b1)$) -- ($1.2*(b1)$);
    \node[orange!85!black, anchor=north west, xshift=-2, yshift=3, rotate=-10] at ($0.5*(b1)$) {$c_r=0$};

    \draw[-{Stealth[length=2.5mm]}, red, very thick] (0,0,0) -- (n);
    
    \def\orthlen{0.15}
    \coordinate (p1) at ($ -\orthlen*(b1) $);
    \coordinate (p2) at ($ \orthlen*(n) $);
    \coordinate (p3) at ($ (p1) + (p2) $);
    \draw[red, thick] (p1) -- (p3) -- (p2);

    \node[text=red, align=left, anchor=east, xshift=10] at (n) 
    {
        Uncontrollable Direction \\
        $n_c = \begin{psmallmatrix} -[v_r]_i \\ 1 \end{psmallmatrix}$
    };
    
    \node[blue, align=center, anchor=south, rotate=-10] at ($-0.2*(b1) - 0.5*(b2)$) {Controllable Plane};

\end{tikzpicture}
    \caption{Input-to-RBM map. Reachable RBMs form the blue ``controllable plane,'' orthogonal to the uncontrollable direction $n_c$. Its intersection with the translation subspace ($c_r=0$, horizontal plane) is the 1D set of perfect-recovery outcomes.}
    \label{fig:unreachable_rbm_illustration}
\end{figure}

\section{Conclusion}
We presented a geometric input-output analysis of hidden modes in distance-based formation control under a single localized actuator. For minimally connected actuators ($\deg(i)=d$), we proved that the uncontrollable subspace is geometrically bounded by the local rotational subspace $\mathcal{T}_i$. {\color{black} This yields an orthogonal decomposition of the uncontrollable subspace whose rigid-body
component is exactly the global rotational subspace $\mathcal{R}_i$, and whose
deformational component is contained in $\mathcal{T}_i \cap \image(R(p^*)^\top)$.} Because all deformational hidden modes in this LTI system are exponentially stable transients, the formation's steady-state response is restricted entirely to the rigid-body eigenspace. This enabled us to prove a shape-recovery dichotomy from a purely local input-output perspective: a localized impulsive disturbance produces persistent shape distortion if and only if its component along the local rotational RBM $[v_r]_i$ is non-zero. Future work will investigate the exact geometric characterization of over-constrained nodes ($\deg(i) > d$), where local equilibrium stresses theoretically may permit hidden modes outside $\mathcal{T}_i$, and extend this input-output framework to the original nonlinear dynamics.



\bibliography{ifacconf}             
  
\end{document}